\documentclass[pra]{revtex4}
 \usepackage{amssymb} \usepackage{graphicx}
\usepackage[all]{xy}

\begin{document}
 \title{Extended superalgebras from twistor and Killing spinors}

\author{\"Umit Ertem}
 \email{umitertemm@gmail.com}
\address{Department of Physics,
Ankara University, Faculty of Sciences, 06100, Tando\u gan-Ankara,
Turkey\\}

\begin{abstract}

The basic first-order differential operators of spin geometry that are Dirac operator and twistor operator are considered. Special types of spinors defined from these operators such as twistor spinors and Killing spinors are discussed. Symmetry operators of massless and massive Dirac equations are introduced and relevant symmetry operators of twistor spinors and Killing spinors are constructed from Killing-Yano (KY) and conformal Killing-Yano (CKY) forms in constant curvature and Einstein manifolds. The squaring map of spinors gives KY and CKY forms for Killing and twistor spinors respectively. They constitute a graded Lie algebra structure in some special cases. By using the graded Lie algebra structure of KY and CKY forms, extended Killing and conformal superalgebras are constructed in constant
curvature and Einstein manifolds.

\end{abstract}

\maketitle

\section{Introduction}

In spin geometry, sections of the spinor bundle on a spin manifold are called spinor fields and basic first-order differential operators defined on spinor fields are Dirac operator and twistor operator. They are defined in terms of the Clifford product and spinor covariant derivative. Spinors that are in the kernel of the Dirac operator are called harmonic spinors and the set of spinors which are in the kernel of the twistor operator consists of twistor spinors \cite{Penrose Rindler, Lichnerowicz1, Baum Leitner}. Moreover, twistor spinors that also correspond to the eigenspinors of the Dirac operator with non-zero eigenvalues are called Killing spinors. The existence of twistor and Killing spinors put some restrictions on the curvature characteristics of the underlying manifold \cite{Baum Friedrich Grunewald Kath, Lichnerowicz2}.

Spinor bilinears that are constructed from a spinor and a dual spinor can be written as a sum of differential forms through Fierz identity. The 1-form part of spinor bilinears which are called Dirac currents correspond to the metric duals of Killing vector fields and conformal Killing vector fields in the cases of Killing spinors and twistor spinors, respectively. Similarly, $p$-form Dirac currents of Killing and twistor spinors satisfy the Killing-Yano (KY) and conformal Killing-Yano (CKY) equations, respectively \cite{Semmelmann, Acik Ertem}. KY forms and CKY forms are defined as the antisymmetric generalizations of Killing and conformal Killing vectors to higher degree differential forms and they satisfy graded Lie algebra structures in constant curvature manifolds \cite{Kastor Ray Traschen, Acik Ertem2,  Ertem1}.

On the other hand, one can construct symmetry operators for the differential equations written in terms of Dirac operator and twistor operator. Symmetry operator of an equation takes a solution of the equation and gives another one. Symmetry operators of massles and massive Dirac equations are constructed from CKY and KY forms, respectively \cite{Benn Charlton, Benn Kress1, Acik Ertem Onder Vercin, Benn Kress2}. On some special cases, these can also be generalized to twistor and Killing spinor cases. In constant curvature manifolds, the symmetry operators for the Killing spinor equation can be constructed from odd degree KY forms and for twistor equation, they are constructed out of CKY forms. Moreover, symmetry operators of the twistor equation can also be constructed from normal CKY forms in Einstein manifolds \cite{Ertem2, Ertem3}.

In this paper, we review the constructions of some superalgebra structures by using solutions of some special spinor equations, their spinor bilinears and symmetry operators.  For extended Killing superalgebras, the odd part of the superalgebra consists of Killing spinors and even part correspond to KY forms. The bilinear operations of the superalgebras are defined as the Lie bracket of KY forms, symmetry operators of Killing spinors and the spinor bilinear map of spinors and they satisfy a superalgebra structure in constant curvature manifolds. They are generalizations of Killing superalgebras which consist of Killing vector fields and Killing spinors. For extended conformal superalgebras, the even and odd parts are CKY forms and twistor spinors, respectively. The bilinear operations correspond to graded Lie bracket of CKY forms, symmetry operators of twistor spinors and spinor bilinear map of twistor spinors. They satisfy a superalgebra structure in constant curvature manifolds and correspond to generalizations of conformal superalgebras which consist of conformal Killing vector fields and twistor spinors. In Einstein manifolds, extended conformal superalgebras can also be defined in terms of normal CKY forms and twistor spinors.

\section{Dirac and twistor operators}

Let us denote the orthonormal coframe basis on a spin manifold $M$ as $\{e^a\}$. They correspond to the metric duals of the frame basis $\{X_a\}$ and we have the relation $e^a(X_b)=\delta^a_b$. In terms of the Clifford product which is denoted by $.$ and the spinor covariant derivative $\nabla_{X_a}$ with respect to any vector field $X_a$, the Dirac operator is defined as follows
\begin{equation}
\displaystyle{\not}D=e^a.\nabla_{X_a}.
\end{equation}
Spinors that are in the kernel of the Dirac operator are called harmonic spinors and satisfy the following massless Dirac equation
\begin{equation}
\displaystyle{\not}D\psi=0
\end{equation}
and the eigenspinors of the Dirac operator with non-zero eigenvalue $m$ satisfy the massive Dirac equation
\begin{equation}
\displaystyle{\not}D\psi=m\psi.
\end{equation}

For any vector field $X$ and its metric dual $\widetilde{X}$, the twistor operator is defined in $n$ dimensions as
\begin{equation}
{\cal{P}}_X=\nabla_X-\frac{1}{n}\widetilde{X}.\displaystyle{\not}D.
\end{equation}
Spinors that are in the kernel of the twistor operator are called twistor spinors and satisfy the following twistor equation
\begin{equation}
\nabla_X\psi=\frac{1}{n}\widetilde{X}.\displaystyle{\not}D\psi.
\end{equation}
The existence of twistor spinors implies some constraints on the curvature characteristics of the ambient manifold. By taking second covariant derivatives of the twistor spinors in (5) and using the definitions of the curvature characteristics as $R_{ab}$ are curvature 2-forms, $P_a=i_{X^b}R_{ba}$ are Ricci 1-forms and ${\cal{R}}=i_{X^a}P_a$ is the curvature scalar, one obtains the following integrability conditions of the twistor equation
\begin{equation}
\nabla_{X_a}\displaystyle{\not}D\psi=\frac{n}{2}K_a.\psi
\end{equation}
\begin{equation}
\displaystyle{\not}D^2\psi=-\frac{n}{4(n-1)}{\cal{R}}\psi
\end{equation}
\begin{equation}
C_{ab}.\psi=0
\end{equation}
where the 1-form $K_a$ is defined as $K_a=\frac{1}{n-2}\left(\frac{\cal{R}}{2(n-1)}e_a-P_a\right)$ and $C_{ab}$ are conformal 2-forms. The manifolds that admits twistor spinors are classified \cite{Lischewski}. Especially, twistor spinors can exist on conformally flat spin manifolds.

Spinors that satisfy both twistor equation and massive Dirac equation are called Killing spinors and from (3) and (5) one can see that they satisfy the following Killing spinor equation
\begin{equation}
\nabla_X\psi=\lambda\widetilde{X}.\psi
\end{equation}
where $\lambda:=m/n$ is the Killing number which is a real or pure imaginary number. The integrability conditions of Killing spinor equation restricts the curvature scalar of the manifold and it reads as
\begin{equation}
{\cal{R}}=-4\lambda^2n(n-1).
\end{equation}
So, Killing spinors can exist on constant curvature or Einstein manifolds and Killing number is real for negative curvature manifolds and pure imaginary for positive curvature manifolds.

\section{Spinor bilinears}

In a spin manifold $M$, let us denote the spinor space as $S$ and the dual spinor space as $S^*$. The tensor product of the spinor space and the dual spinor space gives the algebra of endomorphisms over the spinor space $S\otimes S^*= \text{End} S$ and it corresponds to the Clifford algebra of relevant dimension. So, the spinor bilinears which are the tensor products of a spinor with its dual can be written as a sum of different degree differential forms
\begin{equation}
\psi\otimes\overline{\psi}=(\psi,\psi)+(\psi, e_a.\psi)e^a+(\psi, e_{ba}.\psi)e^{ab}+...+(\psi, e_{a_p...a_2a_1}.\psi)e^{a_1a_2...a_p}+...+(-1)^{\lfloor{n/2}\rfloor}(\psi, z.\psi)z
\end{equation}
where $e^{a_1a_2...a_p}=e^{a_1}\wedge e^{a_2}\wedge...\wedge e^{a_p}$, $\lfloor{ }\rfloor$ is the floor function that takes the integer part of the argument and $z$ is the volume form. $(\,,\,)$ is the spinor inner product an (11) is called as Fierz identitiy. The $p$-form components on the right hand side of (11) are called $p$-form Dirac currents and denoted as
\begin{equation}
(\psi\overline{\psi})_p=(\psi, e_{a_p...a_2a_1}.\psi)e^{a_1a_2...a_p}.
\end{equation}
For $p=1$, the vector fields that are metric duals of them are called Dirac currents.

Dirac currents of Killing spinors correspond to Killing vector fields. Moreover, $p$-form Dirac currents of Killing spinors or their Hodge duals, depending on the degree of the form and the automorphisms of the Clifford algebra that are used in the definition of the spinor inner product, satisfy the KY equation \cite{Acik Ertem}. KY forms are antisymmetric generalizations of Killing vector fields to higher degree differential forms. A KY $p$-form $\omega$ satisfy the following equation
\begin{equation}
\nabla_X\omega=\frac{1}{p+1}i_Xd\omega
\end{equation}
for any vector field $X$ and $d$ denotes the exterior derivative operation, $i_X$ denotes the interior derivative or contruction operation with respect to $X$. So, the $p$-form Dirac current defined in (12) is a squaring map of a spinor which takes a Killing spinor (or two different Killing spinors) and gives a KY form. In constant curvature manifolds, KY forms constitute a graded Lie algebra structure under Schouten-Nijenhuis (SN) bracket which is defined for a $p$-form $\omega_1$ and a $q$-form $\omega_2$ as follows
\begin{equation}
[\omega_1,\omega_2]_{SN}=i_{X^a}\omega_1\wedge\nabla_{X_a}\omega_2+(-1)^{pq}i_{X^a}\omega_2\wedge\nabla_{X_a}\omega_1.
\end{equation}
This gives a $(p+q-1)$-form which satisfies (13) and the graded Lie algebra structure of KY forms arises from the fact that SN bracket has the following properties of (skew)-symmetry and graded Jacobi identity
\begin{eqnarray}
[\omega_1,\omega_2]_{SN}&=&(-1)^{pq}[\omega_2,\omega_1]_{SN}\\
(-1)^{p(r+1)}[\omega_1,[\omega_2,\omega_3]_{SN}]_{SN}&+&(-1)^{q(p+1)}[\omega_2,[\omega_3,\omega_1]_{SN}]_{SN}+(-1)^{r(q+1)}[\omega_3,[\omega_1,\omega_2]_{SN}]_{SN}=0
\end{eqnarray}
where $\omega_3$ is a KY $r$-form \cite{Kastor Ray Traschen, Acik Ertem2}. In particular, one cane see that odd KY forms satisfy a Lie algebra with respect to SN bracket.

On the other hand, Dirac currents of twistor spinors are conformal Killing vector fields and $p$-form Dirac currents of them correspond to CKY forms \cite{Semmelmann, Acik Ertem}. CKY forms are antisymmetric generalizations of conformal Killing vector fields to higher degree differential forms and a CKY $p$-form $\omega$ satisfy the following equation in $n$-dimensions
\begin{equation}
\nabla_X\omega=\frac{1}{p+1}i_Xd\omega-\frac{1}{n-p+1}\widetilde{X}\wedge\delta\omega
\end{equation}
for any vector field $X$ and $\delta$ denotes the co-derivative operation. So, the $p$-form Dirac current defined in (12) is a squaring map of a spinor which takes a twistor spinor (or two different twistor spinors) and gives a CKY form. Indeed, KY forms are a subset of CKY forms which are co-closed $\delta\omega=0$ as can be seen from (13) and (17). In constant curvature manifolds, CKY forms constitute a graded Lie algebra structure with respect to the following bracket which is defined for a CKY $p$-form $\omega_1$ and CKY $q$-form $\omega_2$
\begin{eqnarray}
[\omega_1,\omega_2]_{CKY}&=&\frac{1}{q+1}i_{X^a}\omega_1\wedge i_{X_a}d\omega_2+\frac{(-1)^p}{p+1}i_{X^a}d\omega_1\wedge i_{X_a}\omega_2\nonumber\\
&&+\frac{(-1)^p}{n-q+1}\omega_1\wedge\delta\omega_2+\frac{1}{n-p+1}\delta\omega_1\wedge\omega_2.
\end{eqnarray}
This gives a $(p+q-1)$-form which satisfies (17) and the graded Lie algebra structure of CKY forms arises from the fact that CKY bracket has the following properties of (skew)-symmetry and graded Jacobi identity
\begin{eqnarray}
[\omega_1,\omega_2]_{CKY}&=&(-1)^{pq}[\omega_2,\omega_1]_{CKY}\\
(-1)^{p(r+1)}[\omega_1,[\omega_2,\omega_3]_{CKY}]_{CKY}&+&(-1)^{q(p+1)}[\omega_2,[\omega_3,\omega_1]_{CKY}]_{CKY}+(-1)^{r(q+1)}[\omega_3,[\omega_1,\omega_2]_{CKY}]_{CKY}=0
\end{eqnarray}
where $\omega_3$ is a CKY $r$-form. Moreover, normal CKY forms which are defined by imposing special constraints on the integrability conditions of CKY forms also satisfy a graded Lie algebra structure in Einstein manifolds \cite{Ertem1}.

\section{Symmetry operators}

Symmetry operators are defined as solution generating operators for differential equations. A symmetry operator of an equation takes a solution of the equation and gives another solution. A complete set of mutually commuting symmetry operators of an equation can be used to find general solutions of the equation by the method of separation of variables.

Lie derivative with respect to a Killing vector field $K$ which is denoted by ${\cal{L}}_K$ is a symmetry operator for the massive Dirac equation (3). So, if $\psi$ is a solution of the massive Dirac equation, then ${\cal{L}}_K\psi$ is also a solution. Similarly, ${\cal{L}}_K$ is also a symmetry operator for Killing spinor equation (9). Action of the Lie derivative with respect to $K$ on spinor fields is written as follows \cite{Kosmann}
\begin{equation}
{\cal{L}}_K\psi=\nabla_K\psi+\frac{1}{4}d\widetilde{K}.\psi.
\end{equation}
Moreover, more general first-order symmetry operators of the massive Dirac equation can be written in terms of KY forms \cite{Benn Kress1, Acik Ertem Onder Vercin}. For a KY $p$-form $\omega$, first-order symmetry operators of the massive Dirac equation are defined as
\begin{equation}
L_{\omega}\psi=(i_{X^a}\omega).\nabla_{X_a}\psi+\frac{p}{2(p+1)}d\omega.\psi.
\end{equation}
(22) reduces to (21) for $p=1$. However, different from ${\cal{L}}_K$ case, (22) is a symmetry operator for Killing spinor equation only for odd $p$ and in constant curvature manifolds \cite{Ertem2}. It can also be written in the following form by using the Killing spinor equation (9)
\begin{equation}
L_{\omega}\psi=-(-1)^p\lambda p\omega.\psi+\frac{p}{2(p+1)}d\omega.\psi.
\end{equation}
These can be summarized as follows
\begin{displaymath}
\xymatrix{\displaystyle{\not}D\psi=m\psi \ar[d] & \nabla_X\psi=\lambda\widetilde{X}.\psi \ar[d] \\
{\cal{L}}_K \ar[d]^{\text{all} \,p} & {\cal{L}}_K \ar[d]^{\text{odd} \,p,\,\text{const. curv.}} \\
L_{\omega}=i_{X^a}\omega.\nabla_{X_a}+\frac{p}{2(p+1)}d\omega \ar@/^/[u]^{p=1} & L_{\omega}=i_{X^a}\omega.\nabla_{X_a}+\frac{p}{2(p+1)}d\omega \ar@/^/[u]^{p=1}}
\end{displaymath}

A Killing vector field $K$ preserves the metric structure $g$ on the manifold; ${\cal{L}}_Kg=0$. On the other hand, a conformal Killing vector field $C$ preserves the metric up to a conformal factor; ${\cal{L}}_Cg=2\mu g$ where $\mu$ is a function. Since massless Dirac equation and twistor equation are conformally covariant equations, their symmetry operators are related to the conformal symmetries of the manifold. In $n$-dimensions, the following operator, written in terms of a conformal Killing vector field $C$ and related function $\mu$, is a symmetry operator for the massless Dirac equation
\begin{equation}
{\cal{L}}_C+\frac{1}{2}(n-1)\mu.
\end{equation}
Indeed, more general first-order symmetry operators for the massless Dirac equation can be written in terms of CKY forms \cite{Benn Charlton}. For a CKY $p$-form $\omega$, general first-order symmetry operators are defined as
\begin{equation}
L_{\omega}\psi=(i_{X^a}\omega).\nabla_{X_a}\psi+\frac{p}{2(p+1)}d\omega.\psi-\frac{n-p}{2(n-p+1)}\delta\omega.\psi.
\end{equation}
This reduces to (24) for $p=1$ because of the equality $\delta\widetilde{C}=-n\mu$. Symmetry operators of twistor equation can also be written in terms of conformal Killing vector fields as
\begin{equation}
{\cal{L}}_C-\frac{1}{2}\mu.
\end{equation}
More general first-order symmetry operators of the twistor equation can also be written in terms of CKY forms as in (25). However, in that case we have the following first-order symmetry operator for a CKY $p$-form $\omega$
\begin{equation}
L_{\omega}\psi=(i_{X^a}\omega).\nabla_{X_a}\psi+\frac{p}{2(p+1)}d\omega.\psi+\frac{p}{2(n-p+1)}\delta\omega.\psi
\end{equation}
and the manifold has to have constant curvature. This also reduces to (26) for $p=1$. It is also a symmetry operator of the twistor equation for normal CKY forms in Einstein manifolds \cite{Ertem3}. By using the twistor equation (5), (27) can also be written in the following form
\begin{equation}
L_{\omega}\psi=-(-1)^p\frac{p}{n}\omega.\displaystyle{\not}D\psi+\frac{p}{2(p+1)}d\omega.\psi+\frac{p}{2(n-p+1)}\delta\omega.\psi.
\end{equation}
Conformal cases can be summarized as follows
\begin{displaymath}
\xymatrix{\displaystyle{\not}D\psi=0 \ar[d] & \nabla_X\psi=\frac{1}{n}\widetilde{X}.\displaystyle{\not}D\psi \ar[d] \\
{\cal{L}}_C+\frac{1}{2}(n-1)\mu \ar[d] & {\cal{L}}_C-\frac{1}{2}\mu \ar[d]^{\text{const. curv.},\,\text{Einstein}} \\
L_{\omega}=i_{X^a}\omega.\nabla_{X_a}+\frac{p}{2(p+1)}d\omega-\frac{n-p}{2(n-p+1)}\delta\omega \ar@/^/[u]^{p=1} & L_{\omega}=i_{X^a}\omega.\nabla_{X_a}+\frac{p}{2(p+1)}d\omega+\frac{p}{2(n-p+1)}\delta\omega \ar@/^/[u]^{p=1}}
\end{displaymath}

\section{Extended superalgebras}

The relation between Killing spinors and Killing vector fields gives rise to a superalgebra structure called Killing superalgebra. Similarly, twistor spinors and conformal Killing vector fields also constitute a superalgebra that is conformal superalgebra \cite{Habermann, de Medeiros Hollands, Lischewski}. However, we know that from the squaring map of Killing spinors and twistor spinors higher degree objects such as KY forms and CKY forms arise and this may lead to extensions of the Killing and conformal superalgebras. In this scetion, we show that these extensions are possible in constant curvature and Einstein manifolds through generalized symmetry operators.

A superalgebra $\mathfrak{g}=\mathfrak{g}_0\oplus\mathfrak{g}_1$ is a graded algebra with the even part $\mathfrak{g}_0$ and the odd part $\mathfrak{g}_1$. A bilinear multiplication operation that defines the superalgebra structure is given as follows
\begin{equation}
[.,.]:\mathfrak{g}_i\times\mathfrak{g}_j\longrightarrow\mathfrak{g}_{i+j}
\end{equation}
where $i,j=0,1$ mod 2. For the elements $a,b\in\mathfrak{g}$, the bilinear operation satisfies the following (skew)-symmetry property
\begin{equation}
[a,b]=-(-1)^{|a||b|}[b,a]
\end{equation}
where $|a|$ corresponds to 0 or 1 for $a$ is in $\mathfrak{g}_0$ or $\mathfrak{g}_1$. If the bilinear operation satisfies the graded Jacobi identities, then the superalgebra is called a Lie superalgebra.

Killing superalgebras $\mathfrak{g}=\mathfrak{g}_0\oplus\mathfrak{g}_1$ defined in terms of Killing spinors and Killing vector fields. The even part $\mathfrak{g}_0$ of the superalgebra correspond to the Lie algebra of Killing vector fields and the odd part $\mathfrak{g}_0$ consists of Killing spinors. The bilinear operations of the superalgebra are defined as follows. The even-even bracket which takes two Killing vectors and gives again a Killing vector is the ordinary Lie bracket $[.\,,\,.]$ of vector fields
\begin{equation}
[.\,,\,.]:\mathfrak{g}_0\times\mathfrak{g}_0\longrightarrow\mathfrak{g}_0
\end{equation}
The even-odd bracket which takes a Killing vector and a Killing spinor and gives a Killing spinor is the Lie derivative ${\cal{L}}$ of spinor fields with respect to Killing vector fields since it is a symmetry operator for Killing spinors
\begin{equation}
{\cal{L}}:\mathfrak{g}_0\times\mathfrak{g}_1\longrightarrow\mathfrak{g}_1
\end{equation}
The odd-odd bracket which takes two Killing spinors and gives a Killing vector field is the metric dual of the 1-form part of the squaring map defined in (12) which is the Dirac current of the Killing spinor
\begin{equation}
(\,\,)_1:\mathfrak{g}_1\times\mathfrak{g}_1\longrightarrow\mathfrak{g}_0
\end{equation}
Since these brackets satisfy the (skew)-symmetry condition in (30), we have the Killing superalgebra of Killing vector fields and Killing spinors.

As we have seen in Sections 3 and 4, $p$-form Dirac currents of Killing spinors correspond to KY forms or their Hodge duals and symmetry operators of Killing spinors can be generalized to odd KY forms. This means that we can extend the Killing superalgebra to include KY forms \cite{Ertem2}. Let us denote the extended Killing superalgebra as $\bar{\mathfrak{g}}=\bar{\mathfrak{g}}_0\oplus\bar{\mathfrak{g}}_1$ where the even part $\bar{\mathfrak{g}}_0$ is the Lie algebra of odd KY forms and the odd part $\bar{\mathfrak{g}}_1$ is the space of Killing spinors. Now, we can define the brackets of the extended superalgebra in terms of SN bracket, generalized symmetry operators and $p$-form Dirac currents. Since the graded Lie algebra of KY forms defined in (14) and symmetry operators of Killing spinors written in terms of KY forms in (23), can only be defined in constant curvature manifolds, extended Killing superalgebras can only be constructed in constant curvature manifolds. The even-even bracket corresponds to SN bracket of KY forms since odd KY forms satisfy a Lie algebra with respect to it
\begin{equation}
[.\,,\,.]_{SN}:\bar{\mathfrak{g}}_0\times\bar{\mathfrak{g}}_0\longrightarrow\bar{\mathfrak{g}}_0
\end{equation}
The even-odd bracket is the generalized symmetry operators of Killing spinors defined in (23)
\begin{equation}
L:\bar{\mathfrak{g}}_0\times\bar{\mathfrak{g}}_1\longrightarrow\bar{\mathfrak{g}}_1
\end{equation}
The odd-odd bracket is defined as the $p$-form Dirac currents of Killing spinors written in (12)
\begin{equation}
(\,\,)_p:\bar{\mathfrak{g}}_1\times\bar{\mathfrak{g}}_1\longrightarrow\bar{\mathfrak{g}}_0
\end{equation}
These brackets satisfy the (skew)-symmetry condition in (30) and we can define the extended Killing superalgebra in terms of odd KY forms and Killing spinors in costant curvature manifolds.

We denote conformal superalgebras as $\mathfrak{k}=\mathfrak{k}_0\oplus\mathfrak{k}_1$ which consists of twistor spinors and conformal Killing vectors. The even part $\mathfrak{k}_0$ is the Lie algebra of conformal Killing vector fields and the odd part $\mathfrak{k}_1$ is the space of twistor spinors. The bilinear operations are defined as follows. The even-even bracket is the ordinary Lie bracket of vector fields and conformal Killing vector fields satisfy a Lie algebra with respect to it
\begin{equation}
[.\,,\,.]:\mathfrak{k}_0\times\mathfrak{k}_0\longrightarrow\mathfrak{k}_0
\end{equation}
The even-odd bracket is the symmetry operator of twistor equation which is defined in terms of conformal Killing vectors in (26)
\begin{equation}
{\cal{L}}-\frac{1}{2}\mu:\mathfrak{k}_0\times\mathfrak{k}_1\longrightarrow\mathfrak{k}_1
\end{equation}
The odd-odd bracket is the metric dual of the 1-form part of (12) for twistor spinors, namely the Dirac currents of twistor spinors which correspond to conformal Killing vectors
\begin{equation}
(\,\,)_1:\mathfrak{k}_1\times\mathfrak{k}_1\longrightarrow\mathfrak{k}_0
\end{equation}
These brackets satisfy the (skew)-symmetry condition and they consist the conformal superalgebra of twistor spinors and conformal Killing vector fields.

Similar to the case of Killing superalgebras, we can extend conformal superalgebras to include CKY forms by using the generalized symmetry operators in (28) and the graded Lie algebra structure of CKY forms in (18). Since the graded Lie algebra structure of CKY forms are constructed in constant curvature manifolds for all CKY forms and in Einstein manifolds for normal CKY forms, the extended conformal superalgebras are also defined in those cases \cite{Ertem3}. Let us denote the extended conformal superalgebras as $\bar{\mathfrak{k}}=\bar{\mathfrak{k}}_0\oplus\bar{\mathfrak{k}}_1$ where the even part $\bar{\mathfrak{k}}_0$ consists of the graded Lie algebra of CKY forms and the odd part $\bar{\mathfrak{k}}_1$ is the space of twistor spinors. The bilinear operations of the extended conformal superalgebras are defined as follows. The even-even bracket is the CKY bracket defined in (18)
\begin{equation}
[.\,,\,.]_{CKY}:\bar{\mathfrak{k}}_0\times\bar{\mathfrak{k}}_0\longrightarrow\bar{\mathfrak{k}}_0
\end{equation}
The even-odd bracket is the generalized symmetry operator of twistor equation defined in terms of CKY forms in (28)
\begin{equation}
L:\bar{\mathfrak{k}}_0\times\bar{\mathfrak{k}}_1\longrightarrow\bar{\mathfrak{k}}_1
\end{equation}
The odd-odd bracket is the $p$-form Dirac currents of twistor spinors defined in (12) which correspond to CKY forms
\begin{equation}
(\,\,)_p:\bar{\mathfrak{k}}_1\times\bar{\mathfrak{k}}_1\longrightarrow\bar{\mathfrak{k}}_0
\end{equation}
Since these brackets satisfy the (skew)-symmetry condititon, we construct the extended conformal superalgebras of twistor spinors and CKY forms in constant curvature manifolds or of twistor spinors and normal CKY forms in Einstein manifolds.

\section{Conclusion}

Killing superalgebras and conformal superalgebras play important roles in supergravity theories in mathematical physics \cite{OFarrill HackettJones Moutsopoulos, OFarrill HackettJones Moutsopoulos Simon}. Supergravity Killing spinors which are parallel spinors with respect to supergravity spin connection and define the number of preserved supersymmetries in the theory are related to geometric Killing spinors in Freund-Rubin backgrounds which are solutions of different bosonic supergravity theories that consist of cartesian products of various dimensional anti-de Sitter spaces and spheres. Because of the Killing superalgebras can characterize the geometry of the underlying manifold, they can be used to classify the supergravity backgrounds in various dimensions. Although, the classification problem is solved for maximally supersymmetric backgrounds, the general classification of supergravity backgrounds has not been obtained yet. Similary, conformal superalgebras can characterize the backgrounds for supersymmetric field theories coupled with conformal supergravity \cite{Cassani Klare Martelli Tomasiello Zaffaroni, de Medeiros}.

We construct extended Killing and extended conformal superalgebras by including KY forms and CKY forms in the superalgebra structure. Since KY forms and CKY forms satisfy graded Lie algebra structures in constant curvature manifolds and they can also be used in the construction of symmetry operators of Killing and twistor spinors, the extended superalgebras can be constructed consistently. So, the extended Killing and conformal superalgebras may provide new insights in the classification problem of supergravity and conformal supergravity backgrouns since they include more general geometric structures that characterize the properties of the underlying manifold.



\begin{references}

\bibitem{Acik Ertem} \"{O}. A\c{c}{\i}k, \"{U}. Ertem, Higher-degree Dirac currents of twistor and Killing spinors in supergravity theories, Class. Quantum Grav. \textbf{32} (2015) 175007.

\bibitem{Acik Ertem2} \"{O}. A\c{c}{\i}k, \"{U}. Ertem, Hidden symmetries and Lie algebra structures from geometric and supergravity Killing spinors, Class. Quantum Grav. \textbf{33} (2016) 165002.

\bibitem{Acik Ertem Onder Vercin} \"{O}. A\c{c}{\i}k, \"{U}. Ertem, M. \"{O}nder, A. Ver\c{c}in, First-order symmetries of the Dirac equation in a curved background: a unified dynamical symmetry condition, Class. Quantum Grav. \textbf{26} (2009) 075001.

\bibitem{Baum Friedrich Grunewald Kath} H. Baum, T. Friedrich, R. Grunewald, I. Kath, Twistors and Killing Spinors on Riemannian Manifolds, Teubner, Leipzig, 1991.

\bibitem{Baum Leitner} H. Baum, F. Leitner, The twistor equation in Lorentzian spin geometry, Math. Z. \textbf{247} (2004) 795.

\bibitem{Benn Charlton} I. M. Benn, P. Charlton, Dirac symmetry operators from conformal Killing-Yano tensors, Class. Quantum Grav. \textbf{14} (1997) 1037.

\bibitem{Benn Kress1} I. M. Benn, J. Kress, First-order Dirac symmetry operators, Class. Quantum Grav. \textbf{21} (2004) 427.

\bibitem{Benn Kress2} I. M. Benn, J. Kress, Differential forms relating twistors to Dirac fields, in: Differential Geometry and its Applications, Proceedings of the 10th International Conference DGA 2007, World Scientific Publishing, Singapore, 2008, pp. 573.

\bibitem{Cassani Klare Martelli Tomasiello Zaffaroni} D. Cassani, C. Klare, D. Martelli, A. Tomasiello, A. Zaffaroni, Supersymmetry in Lorentzian curved spaces and holography, Commun. Math. Phys. \textbf{327} (2014) 577.

\bibitem{de Medeiros} P. de Medeiros, Rigid supersymmetry, conformal coupling and twistor spinors, J. High Energy Phys. \textbf{1409} (2014) 032.

\bibitem{de Medeiros Hollands} P. de Medeiros, S. Hollands, Conformal symmetry superalgebras, Class. Quantum Grav. \textbf{30} (2013) 175016.

\bibitem{Ertem2} \"{U}. Ertem, Symmetry operators of Killing spinors and superalgebras in $AdS_5$, J. Math. Phys. \textbf{57} (2016) 042502.

\bibitem{Ertem1} \"{U}. Ertem, Lie algebra of conformal Killing-Yano forms, Class. Quantum Grav. \textbf{33} (2016) 125033.

\bibitem{Ertem3} \"{U}. Ertem, Twistor spinors and extended conformal superalgebras, arXiv:1605.03361 [hep-th] (2016).

\bibitem{OFarrill HackettJones Moutsopoulos} J. Figueroa-O`Farrill, E. Hackett-Jones, G. Moutsopoulos, The Killing superalgebra of 10-dimensional supergravity backgrounds, Class. Quantum Grav. \textbf{24} (2007) 3291.

\bibitem{OFarrill HackettJones Moutsopoulos Simon} J. Figueroa-O`Farrill, E. Hackett-Jones, G. Moutsopoulos, J. Simon, On the maximal superalgebras of supersymmetric backgrounds, Class. Quantum Grav. \textbf{26} (2009) 035016.

\bibitem{Habermann} K. Habermann, The graded algebra and the derivative $L$ of spinor fields related to the twistor equation, J. Geom. Phys. \textbf{18} (1996) 131.

\bibitem{Kastor Ray Traschen} D. Kastor, S. Ray, J. Traschen, Do Killing-Yano tensors form a Lie algebra?, Class. Quantum Grav. \textbf{24} (2007) 3759.

\bibitem{Kosmann} Y. Kosmann, Derivees de Lie des spineurs, Annal. Math. Pura ed Appl. \textbf{91} (1972) 317.

\bibitem{Lichnerowicz2} A. Lichnerowicz, Killing spinors, twistor spinors and Hijazi inequality, J. Geom. Phys. \textbf{5} (1988) 1.

\bibitem{Lichnerowicz1} A. Lichnerowicz, On the twistor spinors, Lett. Math. Phys. \textbf{18} (1989) 333.

\bibitem{Lischewski} A. Lischewski, Conformal superalgebras via tractor calculus, Class. Quantum Grav. \textbf{32} (2014) 015020.

\bibitem{Penrose Rindler} R. Penrose, W. Rindler, Spinors and Space-Time, vol 2, Cambridge University Press, Cambridge, 1984.

\bibitem{Semmelmann} U. Semmelmann, Conformal Killing forms on Riemannian manifolds, Math. Z. \textbf{245} (2003) 503.

\end{references}
\end{document}